\documentclass{article}
\usepackage[latin1]{inputenc}
\usepackage[english]{babel}
\usepackage{amssymb,amsfonts, amsmath, color}
\usepackage[margin=2.7cm]{geometry}
\usepackage[dvipsnames]{xcolor}
\usepackage{graphicx, color, enumerate}
\usepackage[latin1]{inputenc}
\usepackage[active]{srcltx}
\usepackage{mathtools}
\usepackage{tikz}
\usepackage{pgf}
\usepackage{etex}
\usepackage{verbatim}
\usepackage{tikz-3dplot}
\usepackage{pgfkeys}
\usepackage{amsmath}
\usepackage{amsfonts}
\usepackage{amssymb}
\usepackage{amsthm}
\usepackage{algpseudocode}
\usepackage{float}
\usepackage{xcolor}
\usepackage{mathrsfs}
\usepackage{import}
\usepackage{geometry}
\usepackage{fancyhdr}
\usepackage{fp}
\usepackage[colorlinks, citecolor=blue, linkcolor=red]{hyperref}

\newtheorem{theorem}{Theorem}[section]
\newtheorem{proposition}[theorem]{Proposition}

\newtheorem{lemma}[theorem]{Lemma}
\newtheorem{remark}[theorem]{Remark}
\newtheorem{definition}[theorem]{Definition}

\newcommand{\bcl}{\begin{center}}
\newcommand{\ecl}{\end{center}}
\newcommand{\brl}{\begin{right}}
\newcommand{\erl}{\end{right}}
\newcommand{\ben}{\begin{enumerate}}
\newcommand{\een}{\end{enumerate}}
\newcommand{\overliner}{\begin{array}}
\newcommand{\earr}{\end{array}}
\newcommand{\btab}{\begin{tabular}}
\newcommand{\etab}{\end{tabular}}
\newcommand{\bdoc}{\begin{document}}
\newcommand{\edoc}{\end{document}}
\newcommand{\beqy}{\begin{eqnarray}}
\newcommand{\eeqy}{\end{eqnarray}}

\newcommand{\beqi}{\begin{eqnarray*}}
\newcommand{\eeqi}{\end{eqnarray*}}
\newcommand{\bitem}{\begin{itemize}}

\newcommand{\eitem}{\end{itemize}}
\newcommand{\nln}{\newline}
\newcommand{\newt}{\newtheorem}

\newcommand{\pa}{\partial}
\newcommand{\re}{{I\!\!R}}
\newcommand{\Rn}{\R^N}
\newcommand{\xr}{x\in\R }
\newcommand{\x}{\times}
\newcommand{\dyle}{\displaystyle}
\newcommand{\ene}{{I\!\!N}}
\newcommand{\irn}{\int\limits_{\R^N}}
\newcommand{\io}{\int\limits_{\O}}
\newcommand{\meas}{{\rm meas\,}}
\newcommand{\dif}{\nabla_{xy}}
\newcommand{\sign}{{\rm sign}}
\newcommand{\map}{\longrightarrow }
\newcommand{\imp}{\Longrightarrow }
\renewcommand{\div}{\nabla\cdot }
\newcommand{\sen}{{\rm sen\,}}
\newcommand{\tg}{{\rm tg\,}}
\newcommand{\arcsen}{{\rm arcsen\,}}
\newcommand{\arctg}{{\rm arctg\,}}
\newcommand{\supp}{{\textsl supp\ }}
\newcommand{\ity}{\int_{-\iy}^{+\iy}}
\newcommand{\limit}{\lim\limits}
\newcommand{\limi}{\limit_{n\to\infty}}
\newcommand{\sumi}{\sum\limits_{n=1}^{\infty}}
\newcommand{\ulu}{\underline u}
\newcommand{\ulw}{\underline w}
\newcommand{\ulz}{\underline z}
\newcommand{\ulv}{\underline v}
\newcommand{\uls}{\underline s}
\newcommand{\olu}{\overline u}
\newcommand{\olv}{\overline v}
\newcommand{\ols}{\overline s}
\newcommand{\ob}{\overline\b}
\newcommand{\ovar}{\overline\var}
\newcommand{\wv}{\widetilde v}
\newcommand{\wu}{\widetilde u}
\newcommand{\ws}{\widetilde s}
\renewcommand{\a }{\alpha }
\renewcommand{\b }{\beta }
\newcommand{\g }{\gamma}
\newcommand{\G }{\Gamma }
\renewcommand{\d }{\delta }

\newcommand{\D }{\Delta }
\newcommand{\e }{\varepsilon }
\newcommand{\z }{\zeta }
\renewcommand{\l }{\lambda }
\renewcommand{\L }{\Lambda }
\newcommand{\m }{\mu }
\newcommand{\n }{\nabla }
\newcommand{\s }{\sigma }
\newcommand{\Sig }{\Sigma }
\renewcommand{\t }{\tau }
\newcommand{\var }{\varphi }
\renewcommand{\o }{\omega }
\renewcommand{\O }{\Omega }
\newcommand{\R}{{\mathbb{R}}}
\newcommand{\bC}{{\bf C}}
\newcommand{\bZ}{{\bf Z}}
\newcommand{\bN}{{\bf N}}
\newcommand{\bQ}{{\bf Q}}
\newcommand{\bK}{{\bf K}}
\newcommand{\bI}{{\bf I}}
\newcommand{\bv}{{\bf v}}
\newcommand{\bV}{{\bf V}}
\newcommand{\LL}{\mathcal{L}}
\newcommand{\N}{\mathbb{N}}
\DeclareMathOperator{\suppo}{supp} \DeclareMathOperator{\di}{div}

\newenvironment{Proof}{\Rmovelastskip\vskip12pt
plus 1pt \noindent\em\rm}{\hfill {\qed \hskip .2cm}}

\title{Uniqueness for the Schr\"odinger Equation on Graphs with Potential Vanishing at Infinity}

\author{Fabio Punzo\thanks{Dipartimento di Matematica, Politecnico di Milano, Italia (fabio.punzo@polimi.it).}\,\, and Marcello Svagna \thanks{Dipartimento di Matematica, Politecnico di Milano, Italia (marcello.svagna@mail.polimi.it).}}

\begin{document}
\maketitle

\begin{abstract}
  We investigate the uniqueness, in suitable weighted $\ell^p$
  spaces, of solutions to the Schr\"odinger equation with a potential, posed on infinite graphs.  The potential can tend to zero at infinite with a certain rate.
\end{abstract}



\section{Introduction}
The goal of this paper is to prove, under suitable hypotheses, the uniqueness of solutions to the Schr\"odinger-type equation
\begin{equation}\label{SchrodingerEquation}
\Delta u - V u = 0 \quad  \text{ in }  G,
\end{equation}
where $(G,\mu,\omega)$ is an infinite weighted graph with edge weight $\omega$ and node measure $\mu$, $\Delta$ denotes the weighted graph Laplacian and $V$ is a non-negative given function, which is referred to as a potential.

Graphs play a prominent role in various scientific fields due to their ability to model and describe different phenomena. In mathematics, graphs are particularly significant as they embody the concept of relation, consisting of abstract structures in which distinct elements are interconnected. This characteristic underpins their versatility and applicability across a wide range of problems.

Among various applications, the study of PDEs posed on graphs has recently gained momentum. For instance, several studies have focused both on elliptic equations (see, e.g., \cite{AS1, AS2, BiaMePu1, BP3, Ellipt1,Ellipt2, MePu1, MonPuSo, Ellipt8}) and on  parabolic problems (see, e.g., \cite{Parab1, Parab2, GMP, Parab5b, HuangEta, MoPuSo2, Parab5, Parab8}).
Before providing the mathematical framework in which the problem is set, it is worthy to recall some results in literature related to our problem.

\smallskip

Uniqueness of solutions to \textit{Schrodinger-type} equations like (\ref{SchrodingerEquation}) has been deeply studied on Riemannian manifolds (see, e.g. \cite{GrigRiemannManif1, MeRonc, PuManifolds1}; see also \cite{BiaPu}).

Some $\ell^p$-Liouville theorems for the Schr\"odinger equation with $V\equiv 0$ on graphs can be found in \cite{HuaJost, HuaKel1, Masamune1}.
 Furthermore, in \cite{LibroKeller} it is proven that if $u$ solves equation (\ref{SchrodingerEquation}) and $u$ belongs to a weighted $\ell^2$ space with weight function $\phi =  e^{{-\beta d(x,x_0)}}$, with $\beta > 0$, then $u \equiv 0$.

\smallskip

The most similar results to the ours, however, are obtained in \cite{MePu1}, in which the uniqueness of solutions to (\ref{SchrodingerEquation}) in suitable weighted $\ell^p$ spaces is proven. In \cite{MePu1}, weight functions that decay at infinity exponentially are considered. Moreover, a crucial assumption is that the potential function $V$ is bounded away from zero. Now, we remove that hypothesis on $V$. Indeed, we allow the potential $V(x)$ to tend to zero in a controlled way as $x\to \infty$. This fact causes important difficulties.
In our setting, the weight is a negative power of the distance from a fixed reference point. Therefore, with respect to \cite{MePu1}, on the one hand we enlarge the class to which the potential belongs, on the other hand we restrict the class in which the solution is, since the decay at infinity of the weight is slower than the exponential one. However, we can consider also {\it unbounded} solutions (see Remark \ref{oss1}). Let us mention that our assumption on the potential $V$ is optimal (see Remark \ref{oss2}).

The validity of the Liouville property for equation (\ref{SchrodingerEquation}) posed on infinite graphs, with a potential that can tend to zero at infinity has been also addressed in \cite{BiaMePu1, BP3}. However, in  \cite{BiaMePu1} only {\it bounded} solutions are considered, while in \cite{BP3} on the solutions a point-wise growth condition at infinity is imposed, indeed Phragm\`en-Lindel\"of type theorems are established.
\smallskip

 We shall distinguish the cases $p\geq 2$ and $1\leq p<2$. For $p\geq 2$, our proof is based on an appropriate a-priori estimate (see Proposition \ref{proposition 4.3 NUOVA}) and on the construction of a suitable {\em supersolution} to a sort of {\em adjiont} equation (see Lemma \ref{Lemma 5.1 nuova versione}).
 The a-priori estimate is obtained by multiplying the equation by an appropriate test function, and then by integrating by parts once.

 When $1\leq p<2$, we need to make some extra hypotheses on the potential $V$ and on the graph $G$. Indeed, such method works for any $p\geq 1$, but even if $p\geq 2$ we need stronger hypotheses. In this case the a-priori estimate is obtained by integrating by parts twice (see Proposition \ref{second fundamental proposition (4.4)}). As a consequence, we have to consider a different {\em adjoint equation} (see Lemma \ref{New Lemma 6.1}).
The fact that in order to treat also the case $1\leq p<2$ one requires further hypotheses is typical for such type of results; for instance, the same happens on Riemannian manifolds (see \cite{PuManifolds1}).

 \smallskip

  The paper is organized as follows. In Section \ref{mf}, an introductory overview of the graph setting will be provided, together with the basic tools we will rely on in the following parts.
Section \ref{mr} is devoted to the statements of the main results. In Section \ref{mt} the uniqueness result for $p\geq 2$ will be proved. Finally, in Section \ref{alt} we give the proof in the case $1\leq p<2$.

\section{Mathematical framework}\label{mf}\setcounter{equation}{0}
\subsection{Graph Setting}
Let us introduce some basic material concerning graphs (see \cite{LibroGrigo, LibroKeller}). A weighted Graph $(G, \mu, \omega)$ is a triplet, where
\begin{itemize}
    \item $G$ is an arbitrary set, whose elements are called vertices;
    \item $\mu$ : $G \longrightarrow (0, \infty) $ is a function, called \textit{node measure};
    \item $\omega$ : $G\times G \longrightarrow [0,\infty)$ is a function called \textit{edge weight}.
\end{itemize}
In our particular setting $\omega$ is required to satisfy:
\begin{itemize}
    \item (Absence of loops) $\omega(x,x) = 0$ for all $x \in G$;
    \item (Symmetry of the graph) $\omega(x,y) = \omega(x,y)$ for all $(x,y) \in G \times G$;
    \item (Finite sum) $\sum_{y \in G} \omega(x,y) < \infty$ for all $x \in G$.
\end{itemize}
Notice that if $G$ is a finite set then the graph is said to be \textit{finite}, while if $G$ is infinite and countable the graph is called \textit{infinite}.
We write $x \sim y$ if $\omega(x,y) > 0$ and we say that $x$ is \textit{connected} to $y$ or that the two vertices are \textit{adjacent}. A graph is called \textit{locally finite} if, for each vertex $x \in G$ the number of vertices adjacent to $x$ is finite.
Whenever $\omega(x,y) > 0$ the pair $(x,y)$ is called \textit{edge} of the graph. In this case, the two vertices are called \textit{endpoints} of the graph. According to the restrictions reported above, we are dealing with \textit{undirected} graphs (i.e. graphs whose edges do not have orientation).

An important definition is the following one: a sequence $\{x_k\}_{k=0}^n$ of vertices $\in$ G is called a \textit{path} if $x_k \sim x_{k+1}$ $\forall k=0, ..., n-1$.
A graph is then called \textit{connected} if, for any two vertices $x$ and $y \in G$, there exist a path joining $x$ to $y$.\\
It is now possible to define the \textit{degree} and the \textit{weighted  degree} of a vertex $x \in G$ in the following way:
\[
    \operatorname{deg}(x) \coloneqq \sum_{y \in G} \omega(x,y), \hspace{2 cm} \operatorname{Deg}(x) \coloneqq \frac{deg(x)}{\mu(x)}.
\]
\begin{definition}
A \textit{pseudo metric} on the graph $G$ is a map $d: G\times G \longrightarrow [0,\infty)$ such that
\begin{itemize}
    \item $d(x,x) = 0$ \hspace{1.5 cm} for all $x \in G$;
    \item $d(x,y) = d(y,x)$ \hspace{1.5 cm} for all $x,y \in G$;
    \item $d(x,y) \leq d(x,z) + d(y,z)$ \hspace{1.5 cm} for all $x,y,z \in G$.
\end{itemize}
\end{definition}
\hspace{-0.5 cm}Notice that we called $d$ a \textit{pseudo metric} since in some cases it is possible to find vertices $x,y \in G$ s.t. $d(x,y) = 0, \ x \neq y$ (for instance see \cite{LibroKeller}, Example 11.6), so in general it is not necessarily a metric.\\

For any $x_0 \in G$ and $r >0$, we define the ball $B_r(x_0)$ as
\begin{align}  \label{Ball}
     B_r(x_0) := \{x \in G \ s.t. \ d(x, x_0) < r \}.
\end{align}
We want now to introduce now some further definitions regarding the pseudo metric $d$.
\begin{definition}
The \textit{jump size} $s$ of a pseudo metric $d$ is defined as
  \[
    s = \sup \{d(x,y) : x,y \in G, \omega(x,y) > 0 \}.
\]
\end{definition}
\begin{definition}
    \label{Q-INTRINSIC METRIC}
    Let $q \geq 1$, $C_0 > 0$. We say that a pseudo metric d on $(G,\mu,\omega)$ is {\em q-intrinsic} with bound $C_0$ if
\begin{align*}
    \frac{1}{\mu(x)} \sum_{y \in G} \omega(x,y) d^q (x,y) \leq C_0 \hspace{2 cm} \forall x \in G.
\end{align*}
\end{definition}
\hspace{-0.5 cm}If $q = 2, \ C_0 = 1$, then $d$ is usually called \textit{intrinsic}.
\subsection{Difference Operator and Graph Laplacian}
We first define $\mathfrak{F}$ the set of all the functions $f: G \longrightarrow \mathbb{R}$.
\begin{definition} \label{DIFFERENCE OPERATOR}
Given a function $ f \in \mathit{F}$, $\forall x,y \in G$, we define the \textit{difference operator} $\mathit{\nabla_{xy}}$ : $\mathit{F} \longrightarrow \mathbb{R}$ as
\begin{align*}
\mathit{\nabla_{xy}f} := f(y) - f(x).
\end{align*}
\end{definition}
\hspace{-0.5 cm}For any $f,g \in \mathfrak{F}$, the following equality holds
\begin{align*}
    \mathit{\nabla_{xy}(fg)} = f(x) \mathit{\nabla_{xy}(g)}+ g(y)\mathit{\nabla_{xy}(f)}.
\end{align*}
This formula is commonly known as the \textit{product rule}.
\begin{definition}
\label{Weigthed Graph Laplacian Definition}
Let $(G,\mu,\omega)$ a connected weighted graph, for any function $f \in \mathfrak{F}$ we define
    \begin{align*}
        \Delta f(x) \coloneqq  \frac{1}{\mu(x)}  \sum_{y \sim x} \omega(x,y) [f(y) - f(x)]\quad \forall x\in G.
    \end{align*}
\end{definition}

Let $(G,\mu,\omega)$ be a \textit{locally finite} and connected weighted graph, let $\Omega$ be a finite subset of $G$, then, $\forall f,g \ \in \mathfrak{F}$, the following identity holds
\begin{equation*}
    \begin{aligned}
        &\sum_{x \in \Omega} \Delta f(x)  g(x) \mu(x) = -\frac{1}{2} \sum_{x,y \in \Omega} \nabla_{xy}f\,  \nabla_{xy}g\, \omega(x,y) \\
        &+ \sum_{x \in \Omega} \sum_{y \in \Omega^{c}} \nabla_{xy}f\, g(x) \, \omega(x,y).
    \end{aligned}
\end{equation*}
This identity is commonly known as the Green's formula on (weighted) graphs.

Actually, if we consider at least one of the two functions $f$ and $g$ with finite support coincident with $\Omega$, the Green's formula can be written as
\begin{align}\label{Green's Formula}
    \sum_{x \in \Omega} \Delta f(x) g(x) \mu(x) = -\frac{1}{2} \sum_{x,y \in \Omega} \nabla_{xy}f \, \nabla_{xy}g\,  \omega(x,y).
\end{align}

In the following two Lemmas we recall two properties of the Laplacian (see, e.g., \cite{MePu1}).

\begin{lemma} {}\label{Lemma Prodotto due funz}
Let $f,g \in \mathfrak{F}$. Then, $\forall x \in G$ , one has that
\begin{align*}
    \Delta[f(x)g(x)] = f(x)\Delta g(x) + g(x)\Delta f(x) + \frac{1}{\mu(x)}\sum_{y \in G} (\nabla_{xy} f) (\nabla_{xy} g) \omega(x,y).
\end{align*}
\end{lemma}

\begin{lemma} {}\label{Laplacian of a convex function}
    Let $\psi \in C^1(\mathbb{R})$ be a convex function, and let $u \in \mathfrak{F}$. Then it holds that
    \begin{align*}
        \Delta \psi(u(x)) \geq \psi'(u(x))\Delta u(x) \hspace{1.5 cm} \forall x \in G.
    \end{align*}
\end{lemma}

\subsection{Weighted $\ell^p$ spaces}
\begin{definition}
    \label{Weighted_spaces_definition}
   Let $ \zeta \ : G \longrightarrow (0,\infty) $ be a positive function, $\forall p  \in [1,\infty]$ we define the weighted $\ell_\zeta^p$ space as
   \begin{align*}
       \ell_\zeta^p (G,\mu) = \big\{ u \in \mathfrak{F}  :  \sum_{x \in G} |u(x)|^p \zeta(x) \mu(x) < \infty \big\}, \hspace{1 cm} \text{ for }  p \in [1,\infty),
   \end{align*}
   and
   \begin{align*}
       \ell_\zeta^\infty (G,\mu) = \big\{  u \in \mathfrak{F} \ : \sup_{x \in G} |u(x)| < \infty \big\}.
   \end{align*}
   Notice that, when $\zeta$ is identically $\equiv$ 1 on G, then the above definition coincides with the one of the classical $\ell^p$ space on graphs.
\end{definition}

\hspace{-1 cm}

Given a vertex $x_0 \in G$, we define the weight function $\zeta_\beta$ as
\begin{align}\label{ZETA FUNZ PESO NUOVA}
    \zeta_\beta(x) = [d(x,x_0)+k]^{-\beta} \hspace{1.5 cm} \forall x \in G,
\end{align}
where both $\beta$ and $k$ are two positive parameters.

\section{Main results}\label{mr}\setcounter{equation}{0}
Let us briefly recall the fundamental hypotheses we will always assume from now on
\vspace{1 ex}
\begin{equation}
\label{GENERAL HYPOTHESIS}
\begin{aligned}[t]
& \begin{split}
&\text{ \textit{(i)} $(G,\mu,\omega) $ is an \textit{infinite}, \textit{weighted} and \textit{connected} graph;}
\end{split} \\
& \begin{split}
&\text{\textit{(ii)} there exists a pseudo-metric $d$ with finite jump size $s$;}
\end{split} \\
& \begin{split}
&\text{\textit{(iii)} the ball $B_r(x)$ is a finite set, $\forall x \in G$, $r>0$}.
\end{split} \\
\end{aligned}
\vspace{1 ex}
\end{equation}
Let us underline firstly that hypothesis \textit{(ii)} and \textit{(iii)} imply that the graph is locally finite.

\text{Given $x_0 \in G$, we assume that the function $V \in \mathfrak{F}$ fulfills}
\begin{align}\label{NEW Potential}
    V(x) \geq  c_0[d(x,x_0) + k]^{-\alpha} \hspace{1.5 cm} \forall x \in G,
\end{align}
for some $c_0 > 0$, $k > s$ and $\alpha > 0$.\\

We can now establish our uniqueness result in $\ell_{\zeta_\beta}^p$ with $p\geq 2$.

\begin{theorem} {} \label{new teorema 1}
Assume that $d$ is an intrinsic pseudo metric with related finite jump $s$. Let $u$ be a solution of (\ref{SchrodingerEquation}) with $V$ fulfilling (\ref{NEW Potential}) for some $\alpha \leq 2$.
Suppose that $u \in \ell_{\zeta_\beta} ^p(G,\mu)$ for some $\beta>0$ and $p\geq 2$ fulfilling
\begin{equation}\label{e1f}
   p\geq \frac{\beta^2}{2 c_0} \bigg(\frac{k^{\beta + 1}}{(k-s)^{\beta + 1}}\bigg)^2\,.
\end{equation}
Then $u=0$ in  $G$.
\end{theorem}

\begin{remark}\label{oss1}
Suppose that for some $\bar C>0, x_0\in G, m\in\mathbb N$
\begin{equation}\label{e6f}
\mu(B_R(x_0):=\sum_{x\in B_R(x_0)} \mu(x)\leq \bar C R^{m}\quad \text{ for all } R\geq 2\,.
\end{equation}
In addition assume that, for some $\hat C>0, \sigma>0$,
\[u(x) \leq \hat C [d(x, x_0)+1]^{\sigma} \quad \text{ for all } x\in G.\]
Then it is direct to see that $u\in \ell_{\zeta_\beta} ^p(G,\mu)$, whenever
\[\beta>m+1, \quad \sigma<\frac{\beta-m-1}{p}\,.\]
Therefore, also {\em unbounded} solutions can be considered in Theorem \ref{new teorema 1}.

\smallskip
Note that, obviously, \eqref{e6f} is fulfilled when $G$ is the standard integer lattice $\mathbb Z^m$.
\end{remark}

\begin{remark}\label{oss2}
Observe that hypothesis  (\ref{NEW Potential}) with $\alpha\leq 2$ is optimal. In fact, in \cite[Theorem 3.12]{BP3} it is shown that if $G=\mathbb Z^n$ and
$$V(x)\leq [1+ d(x, x_0)]^{-\alpha}\quad \text{ for all } x\in G,$$ for some $\alpha>2$, then equation \eqref{SchrodingerEquation} admits {\em infinitely many} bounded solutions.
\end{remark}

Now we state our uniqueness result in $\ell_{\zeta_\beta}^p$ with $p\geq 1$. In this case we need some extra assumptions on the potential $V$ and on the pseudo distance $d$ on the graph $G$.

\begin{theorem} {} \label{teo2}
Assume that $d$ is an intrinsic pseudo metric with related finite jump $s$, which is also $1-$intrnsic with bound $C_0$. Let $u$ be a solution of (\ref{SchrodingerEquation}) with $V$ fulfilling (\ref{NEW Potential}) for some $\alpha \leq 1$.
Suppose that $u \in \ell_{\zeta_\beta} ^p(G,\mu)$ for some $\beta>0$ and $p\geq 1$ fulfilling
\begin{equation}\label{e2f}
    p> \frac{C_0}{c_0} \frac{k^{\beta+1}}{(k-s)^{\beta + 1}}\beta.
\end{equation}
Then $u=0$ in  $G$.
\end{theorem}

\begin{remark}
Let us mention that obviously when $\beta=0, 1\leq p_1<p_2$, then $\ell_{\zeta_\beta} ^{p_1}(G,\mu)\subset \ell_{\zeta_\beta} ^{p_2}(G,\mu)$. However, that inclusion fails if $\beta>0$. In that case
\[\ell_{\zeta_{\beta_1}} ^{p_1}(G,\mu)\subset \ell_{\zeta_{\beta_2}} ^{p_2}(G,\mu).\]
provided that $\displaystyle 1\leq p_1<p_2, \beta_2\geq \beta_1\frac{p_2}{p_1}$. For this reason, we had to consider Theorem \ref{teo2} separately from Theorem \ref{new teorema 1}.
\end{remark}

\section{Proof of Theorem \ref{new teorema 1}}\label{mt}

We would like to provide in advance an estimate which will be
used in the following.
\begin{lemma}
\label{Lemma Differenza D(x) e D(y)}
Let $(G,\omega,\mu)$ be an infinite weighted graph, with a pseudo metric $d$ and a finite jump size $s$. Given any $x,y \in G$, with $\omega(x,y) > 0$, one has that
\begin{align}
    \big | [d(y,x_0) + k]^{-\beta} - [d(x,x_0) + k]^{-\beta} \big | \leq \beta d(x,y) [d(x,x_0) + k - s]^{-\beta -1}.
\end{align}
\end{lemma}

\begin{proof}
Fix any $x,y \in G  \ \text{with} \ \omega(x,y) >0$.
We set $t:=d(x,x_0) + k$  and $s:=d(y,x_0) + k$, with $s,t \in \mathbb{R} $ and $s,t > 0$. Then we set $f(t) = t^{-\beta} $ and $f(s) = s^{-\beta} $.
It is then possible to apply the Lagrange Theorem on $f$, which immediately gives
\begin{align*}
    f(s)-f(t) = f'(\xi)(s-t) \hspace{1.5 cm} \text{for some $\xi \in \mathbb{R}$ between $s$ and $t$}.
\end{align*}
Thus
\begin{align*}
    [d(y,x_0) + k]^{-\beta} - [d(x,x_0) + k]^{-\beta} = - \beta \xi^{-\beta -1} \bigg\{[d(y,x_0) + k]-[d(x,x_0) + k] \bigg\},
\end{align*}
where $\xi$ is an intermediate value between $[d(y,x_0) + k]$ and $[d(x,x_0) + k]$. \\
Therefore
\begin{align*}
  \bigg | [d(y,x_0) + k]^{-\beta} - [d(x,x_0) + k]^{-\beta} \bigg | = \bigg | \beta \xi^{-\beta -1} \bigg\{[d(y,x_0) + k]-[d(x,x_0) + k] \bigg\} \bigg |.
\end{align*}
Exploiting the triangular inequality, we obtain
\begin{align}\label{intermediate step of dimostration of distance diff}
    \bigg | [d(y,x_0) + k]^{-\beta} - [d(x,x_0) + k]^{-\beta} \bigg | \leq  \beta \xi^{-\beta -1} d(x,y).
\end{align}
Now, we have to find a good estimate for $\xi^{-\beta-1}$. In order to do this, we can rely on the fact that
\begin{equation*}
    \begin{aligned}
    &\text{if $d(x,x_0) > d(y,x_0)$, then} \hspace{1 cm} [d(y,x_0) + k]^{-\beta-1} > \xi^{-\beta-1};\\
    &\text{if $d(y,x_0) > d(x,x_0)$, then} \hspace{1 cm} [d(x,x_0) + k]^{-\beta-1} > \xi^{-\beta-1}.
    \end{aligned}
\end{equation*}
Now, in the case in which $d(x,x_0) > d(y,x_0)$, exploiting the triangular inequality and the definition of jump size $s$, one has that
\begin{align*}
    d(y,x_0) \geq d(x,x_0) - d(x,y) \hspace{0.5 cm} \text{and consequently } \hspace{0.5 cm}  d(y,x_0) \geq d(x,x_0) - s.
\end{align*}
Recalling then the previous expression, it is then possible to write
\begin{equation*}
    \begin{aligned}
    &\text{if $d(x,x_0) > d(y,x_0)$, then} \hspace{1 cm} [d(x,x_0) + k-s]^{-\beta-1} > \xi^{-\beta-1}:\\
    &\text{if $d(y,x_0) > d(x,x_0)$, then} \hspace{1 cm} [d(x,x_0) + k]^{-\beta-1} > \xi^{-\beta-1}.
    \end{aligned}
\end{equation*}
Since $[d(x,x_0) + k-s]^{-\beta-1} > [d(x,x_0) + k]^{-\beta-1} $ for any $x \in G$, we have that, independently from the relation between $d(x,x_0)$ and $d(y,x_0)$,
\begin{align*}
    [d(x,x_0) + k-s]^{-\beta-1} > \xi^{-\beta-1}.
\end{align*}
It is now possible to recover (\ref{intermediate step of dimostration of distance diff}) and to write
\begin{align*}
    \bigg | [d(y,x_0) + k]^{-\beta} - [d(x,x_0) + k]^{-\beta} \bigg | \leq  \beta [d(x,x_0) + k-s]^{-\beta-1} d(x,y).
\end{align*}
\end{proof}

In the following proposition we establish an important a priori estimate for the solutions of (\ref{SchrodingerEquation}).

\begin{proposition}{}\label{proposition 4.3 NUOVA}
    Let $u$ be a solution to (\ref{SchrodingerEquation}), and let $\eta$ and $\zeta_\beta$ $\in \mathfrak{F}$.
    Suppose moreover that
    \begin{itemize}
        \item $\eta \geq 0$ with finite support;
        \item $[\eta^2(y) -\eta^2(x)][\zeta_\beta(y)-\zeta_\beta(x)] \geq 0$ \ $\forall x,y \in G$ s.t. $x \sim y$.
    \end{itemize}
    Then, for any $p \geq 2$, the following inequality holds
    \begin{equation}\label{e3f}
    \begin{aligned}
        & \frac{1}{2}\sum_{x \in G} |u(x)|^p \eta^2(x) \zeta_\beta(x) \bigg\{ V(x)p\mu(x) - \frac{1}{2} \sum_{y \in G} \omega(x,y) \bigg[ 1-\frac{\zeta_\beta(y)}{\zeta_\beta(x)} \bigg]^2 \bigg\} \\
        & \leq \sum_{x,y \in G} |u(x)|^p \zeta_\beta(y) [\eta(y)-\eta(x)]^2\omega(x,y).
    \end{aligned}
    \end{equation}
\end{proposition}

\begin{proof}
The thesis follows directly from \cite[Proposition 4.3]{MePu1}. In fact, in \cite[Proposition 4.3]{MePu1} it is shown that if
$\eta, \xi$ $\in \mathfrak{F}$,

    \begin{itemize}
        \item $\eta \geq 0$ with finite support;
        \item $[\eta^2(y) -\eta^2(x)][e^{\xi(y)}-e^{\xi(x)}] \geq 0$ \ $\forall x,y \in G$ s.t. $x \sim y$,
    \end{itemize}
    then, for any $p \geq 2$,
    \begin{equation}\label{e4f}
    \begin{aligned}
        & \frac{1}{2}\sum_{x \in G} |u(x)|^p \eta^2(x) e^{\xi(x)} \bigg\{ V(x)p\mu(x) - \frac{1}{2} \sum_{y \in G} \omega(x,y) \bigg[ 1-e^{\xi(y)-\xi(x)}\bigg]^2 \bigg\} \\
        & \leq \sum_{x,y \in G} |u(x)|^p e^{\xi(y)} [\eta(y)-\eta(x)]^2\omega(x,y).
    \end{aligned}
    \end{equation}
By choosing $\xi(x)=\log\zeta_\beta(x)$ in \eqref{e4f}, we derive \eqref{e3f}.
\end{proof}

Now, before starting with the proof of the theorem, we need to state and prove a useful lemma.
\begin{lemma} {}\label{Lemma 5.1 nuova versione}
Let $p \geq 2$. Then
\begin{equation}\label{e5f}
\begin{aligned}
    &\frac{1}{2}\sum_{y \in G}\omega(x,y) \bigg[ 1- \frac{\zeta_\beta(y)}{\zeta_\beta(x)}\bigg]^2 - p V(x) \mu(x)\\
    &\leq \mu(x) \bigg( \frac{1}{2} C_1^2 \beta^2 \big[d(x,x_0) + k]^{-2}  - p c_0[d(x,x_0) + k]^{-\alpha} \bigg).
\end{aligned}
\end{equation}
for all $x \in G$, where
\begin{align}\label{C1 DECLARE}
    C_1 \coloneqq \frac{k^{\beta + 1}}{(k-s)^{\beta + 1}}.
\end{align}

\end{lemma}
\begin{proof}
With a simple grouping of terms, in view of the very definition of $\zeta_\beta$, it is possible to write
\[
\begin{aligned}
     &\frac{1}{2}\sum_{y \in G}\omega(x,y) \bigg[ 1- \frac{\zeta_\beta(y)}{\zeta_\beta(x)}\bigg]^2\\
     & = \frac{1}{2}\sum_{y \in G}\omega(x,y) \bigg[\frac{1}{(d(x,x_0)+k)^{-\beta}} \big((d(x,x_0)+k)^{-\beta} - (d(y,x_0)+k)^{-\beta} \big) \bigg]^2.
\end{aligned}
\]
Consequently,
\[
\begin{aligned}
    &\frac{1}{2}\sum_{y \in G}\omega(x,y) \bigg[ 1- \frac{(d(y,x_0)+k)^{-\beta}}{(d(x,x_0)+k)^{-\beta}}\bigg]^2\\
    & = \frac{1}{2}\sum_{y \in G}\omega(x,y) \frac{1}{(d(x,x_0)+k)^{-2\beta}} \big[(d(x,x_0)+k)^{-\beta} - (d(y,x_0)+k)^{-\beta} \big]^2.
\end{aligned}
\]
From Lemma \ref{Lemma Differenza D(x) e D(y)}, we deduce that
\[
\begin{aligned}
    &\frac{1}{2}\sum_{y \in G}\omega(x,y) \bigg[ 1- \frac{(d(y,x_0)+k)^{-\beta}}{(d(x,x_0)+k)^{-\beta}}\bigg]^2\\
    & \leq \frac{1}{2}\sum_{y \in G}\omega(x,y) \frac{1}{(d(x,x_0)+k)^{-2\beta}} \big[\beta d(x,y) [d(x,x_0) + k - s]^{-\beta -1}]^2.
\end{aligned}
\]

\smallskip
Now, we want to find a constant $C_1$ such that
\begin{align} \label{costante C1 origine}
    [d(x,x_0) + k-s]^{-\beta-1} \leq C_1 [d(x,x_0) + k]^{-\beta-1} \hspace{1 cm} \forall x  \in G.
\end{align}
In order to do this, let us introduce the function $g : \mathbb{R^+} \longrightarrow \mathbb{R}$,
\begin{align*}
    g(t) \coloneqq \frac{[t+k-s]^{-\beta-1}}{[t+k]^{-\beta-1}}, \quad t\geq 0.
\end{align*}
Thus
\begin{equation*}
\begin{aligned}
    & g'(t) = \frac{(\beta+1) [t + k]^{-\beta-2} [t + k-s]^{-\beta-1} - (\beta+1)[t + k-s]^{-\beta-2}[t + k]^{-\beta-1}}{[t + k]^{-2(\beta +1)}}.
\end{aligned}
\end{equation*}
Grouping the terms one has
\begin{align*}
    \frac{(\beta+1)[t + k]^{-\beta-2}[t + k-s]^{-\beta-2} \bigg[t + k-s - t - k \bigg]}{[t + k]^{-2(\beta +1)}}.
\end{align*}
We see immediately that the numerator is a negative term for every value of $t$ in the domain, while the denominator is instead always positive. The function $g(t)$ is then decreasing,
this means that when $t = 0$ this function reaches its maximum value. We can apply this analysis done on the function $g$ to the ratio $\frac{ [d(x,x_0) + k-s]^{-\beta-1}}{[d(x,x_0) + k]^{-\beta-1}}$, thus saying that, for $d(x,x_0) = 0$, this expression reaches its maximum value.
Consequently, the optimal choice for $C_1$ is that in \eqref{C1 DECLARE}.

\smallskip

By virtue of (\ref{costante C1 origine}) and (\ref{C1 DECLARE}), we can infer that
\begin{equation*}
\begin{aligned}
    &\frac{1}{2}\sum_{y \in G}\omega(x,y) \bigg[ 1- \frac{(d(y,x_0)+k)^{-\beta}}{(d(x,x_0)+k)^{-\beta}}\bigg]^2\\
    & \leq \frac{1}{2}\sum_{y \in G}\omega(x,y) \frac{1}{(d(x,x_0)+k)^{-2\beta}} \big[\beta d(x,y) C_1 [d(x,x_0) + k]^{-\beta-1}]^2.
\end{aligned}
\end{equation*}
Now, simplifying the proper terms and bringing outside the sum the constants, one can obtain
\begin{equation*}
\begin{aligned}
    &\frac{1}{2}\sum_{y \in G}\omega(x,y) \bigg[ 1- \frac{(d(y,x_0)+k)^{-\beta}}{(d(x,x_0)+k)^{-\beta}}\bigg]^2\\
    & \leq \frac{1}{2} C_1^2 \beta^2 \big[(d(x,x_0) + k)]^{-2}\sum_{y \in G}\omega(x,y) d^2(x,y).
\end{aligned}
\end{equation*}
Exploiting the intrinsic property of $d$, we have that
\begin{equation*}
    \frac{1}{2}\sum_{y \in G}\omega(x,y) \bigg[ 1- \frac{(d(y,x_0)+k)^{-\beta}}{(d(x,x_0)+k)^{-\beta}}\bigg]^2 \leq \frac{1}{2} C_1^2 \beta^2 \big[(d(x,x_0) + k)]^{-2} \mu(x).
\end{equation*}
It is then possible to write that
\begin{equation*}
    \begin{aligned}
        &  \frac{1}{2}\sum_{y \in G}\omega(x,y) \bigg[ 1- \frac{(d(y,x_0)+k)^{-\beta}}{(d(x,x_0)+k)^{-\beta}}\bigg]^2 - p V(x) \mu(x) \\
        & \leq \frac{1}{2} C_1^2 \beta^2 \big[(d(x,x_0) + k)]^{-2} \mu(x) - p V(x) \mu(x).
    \end{aligned}
\end{equation*}
We can now exploit (\ref{NEW Potential}), thus obtaining
\begin{equation*}
    \begin{aligned}
        &  \frac{1}{2}\sum_{y \in G}\omega(x,y) \bigg[ 1- \frac{(d(y,x_0)+k)^{-\beta}}{(d(x,x_0)+k)^{-\beta}}\bigg]^2 - p V(x) \mu(x) \\
        & \leq \frac{1}{2} C_1^2 \beta^2 \big[(d(x,x_0) + k)]^{-2} \mu(x) - p c_0[d(x,x_0) + k]^{-\alpha} \mu(x).
    \end{aligned}
\end{equation*}
And then, grouping the term $\mu(x)$, we arrive to \eqref{e5f}.
\end{proof}

Let us now introduce a \textit{cut-off} function, which will be employed in the proof of the theorem. \\
Let $0<\delta<1$ and $R >0$ , then
\begin{align} \label{definition of eta}
    \eta(x) \coloneqq \min \bigg\{ \frac{[R-s-d(x_0,x)]_+}{\delta R} , 1 \bigg\} \hspace{1.5 cm} \text{for any } x \in G.
\end{align}
Now, by minor changes in the proof of  \cite[Lemma 2.3]{HuangEta}, it is possible to obtain the following result.
\begin{lemma} {}\label{Lemma 5.2 articolo parte due}
    The function $\eta$ defined as in (\ref{definition of eta}) satisfies the  inequality
    \begin{align*}
    |\nabla_{xy}\eta| \leq \frac{1}{\delta R}d(x,y)\chi_{\{(1-\delta)R-2s<d(x,x_0)\leq R \} },
    \end{align*}
    for any $x, y \in G$, $\omega(x,y) >0$. And, consequently, it also satisfies
    \begin{align*}
     |\Delta \eta(x)| \leq \frac{C_0}{\delta R}\chi_{\{(1-\delta)R-2s<d(x,x_0)\leq R \} },
    \end{align*}
    for any $x \in G$.
\end{lemma}

\begin{proof}[Proof of Theorem \ref{new teorema 1}]  It is easy to see that $\eta$ and $\zeta_\beta$ satisfy the hypothesis of the Proposition \ref{proposition 4.3 NUOVA}. Indeed, according to their definition, both  are positive functions, decreasing w.r.t. $d(.\ ,x_0)$. It is then straightforward that if $d(y,x_0) \geq d(x,x_0) $ then both $\eta(y) \leq \eta(x)$ and $\zeta_\beta(y) \leq \zeta_\beta(x)$ and vice-versa. Then, the relation $[\eta^2(y) -\eta^2(x)][\zeta_\beta(y)-\zeta_\beta(x)] \geq 0$ is always true. Therefore, we get
\begin{equation}\label{espressione finale prop 4.3 nuova per dim}
    \begin{aligned}
        & \frac{1}{2}\sum_{x \in G} |u(x)|^p \eta^2(x) \zeta_\beta(x) \bigg\{ V(x)p\mu(x) - \frac{1}{2} \sum_{y \in G} \omega(x,y) \bigg[ 1-\frac{\zeta_\beta(x)}{\zeta_\beta(y)} \bigg]^2 \bigg\} \\
        & \leq \sum_{x,y \in G} |u(x)|^p \zeta_\beta(y) [\eta(y)-\eta(x)]^2\omega(x,y).
    \end{aligned}
\end{equation}
Take
\begin{align*}
    R > \max \bigg\{ \frac{2s}{1-2\delta},1 \bigg\}.
\end{align*}
Observe that
\begin{align}\label{eta=1 ci serve}
    \eta(x) = 1 \quad \text{ whenever } d(x,x_0) \leq \delta R.
\end{align}

By virtue of Lemma \ref{Lemma 5.1 nuova versione}, we have that
\begin{equation}\label{risultato 6.1 con W}
    \begin{aligned}
        &  \frac{1}{2}\sum_{y \in G}\omega(x,y) \bigg[ 1- \frac{(d(y,x_0)+k)^{-\beta}}{(d(x,x_0)+k)^{-\beta}}\bigg]^2 - p V(x) \mu(x) \\
        & \leq \mu(x)\bigg( \frac{1}{2} C_1^2 \beta^2 \big[(d(x,x_0) + k)]^{-2}  - p c_0[d(x,x_0) + k]^{-\alpha} \bigg).
    \end{aligned}
\end{equation}
Let us set $$W:= \sup_{x \in G} \bigg( \frac{1}{2} C_1^2 \beta^2 \big[(d(x,x_0) + k)]^{-2}  - p c_0[d(x,x_0) + k]^{-\alpha} \bigg).$$ Since $\alpha\leq 2,$ from \eqref{e1f} it follows that $W<0$.
Recalling the expression (\ref{espressione finale prop 4.3 nuova per dim}), and exploiting (\ref{eta=1 ci serve}) and (\ref{risultato 6.1 con W}) we are able to find a lower bound for the l.h.s. of (\ref{espressione finale prop 4.3 nuova per dim}) thus obtaining
\begin{equation}\label{lower bound proof nuovo 1 thm}
    \begin{aligned}
        & \frac{1}{2}\sum_{x \in G} |u(x)|^p \eta^2(x) \zeta_\beta(x) \bigg\{ V(x)p\mu(x) - \frac{1}{2} \sum_{y \in G} \omega(x,y) \bigg[ 1-\frac{\zeta_\beta(x)}{\zeta_\beta(y)} \bigg]^2 \bigg\} \\
        & \geq |W| \sum_{x \in B_{\delta_R}} |u(x)|^p \mu(x).
    \end{aligned}
\end{equation}

What we need now, in order to complete the proof, is to find an upper bound to the r.h.s. of (\ref{espressione finale prop 4.3 nuova per dim}) and then to combine the result obtained with the one stated above.

From (\ref{ZETA FUNZ PESO NUOVA}), one has that $\zeta_\beta(y)$ = $[d(y,x_0) + k]^{-\beta}$.
Using a similar reasoning to the one used in the proof of Lemma \ref{Lemma Differenza D(x) e D(y)} we know that, independently from the relation between $d(x,x_0)$ and $d(y,x_0)$, $[d(y,x_0) + k]^{-\beta} \leq [d(x,x_0) + k-s]^{-\beta}$, if $x \sim y$. It is then possible to write
\begin{equation}\label{first step for upper bound nuovo 1 thm}
\begin{aligned}
    &\sum_{x,y \in G} |u(x)|^p \zeta_\beta(y) [\eta(y)-\eta(x)]^2\omega(x,y) \\
    &\leq \sum_{x,y \in G} |u(x)|^p [d(x,x_0) + k-s]^{-\beta}  [\eta(y)-\eta(x)]^2\omega(x,y).
\end{aligned}
\end{equation}
Obviously, for some $C_2>0$, 
\begin{align*}
    [d(x,x_0) + k-s]^{-\beta} \leq C_2 [d(x,x_0) + k]^{-\beta} \quad \text{ for any } x\in G.
\end{align*}
We proceed by recalling Lemma \ref{Lemma 5.2 articolo parte due}, which allows us again to use an useful estimate for the term $[\eta(y)-\eta(x)]^2$, thus allowing us to write
\begin{equation*}
\begin{aligned}
    &\sum_{x,y \in G} |u(x)|^p \zeta_\beta(y) [\eta(y)-\eta(x)]^2\omega(x,y) \\
    &\leq \frac{C_2}{\delta^2 R^2}\sum_{x \in G} \mu(x)\chi_{\{(1-\delta)R-2s \leq d(x,x_0)\leq R \} } |u(x)|^p [d(x,x_0) + k]^{-\beta}.
\end{aligned}
\end{equation*}
Due to the presence of the function $\chi$ inside the sum, the following inequality holds
    \begin{equation}\label{upperbound proof nuovo primo teorema}
    \begin{aligned}
     & \frac{C_2}{\delta^2 R^2}\sum_{x \in G} \mu(x)\chi_{\{(1-\delta)R-2s \leq d(x,x_0)\leq R \} } |u(x)|^p [d(x,x_0) + k]^{-\beta}\\
     &\leq \frac{C_2}{\delta^2 R^2}\sum_{x \in B_{R(x_0)}} \mu(x) |u(x)|^p [d(x,x_0) + k]^{-\beta}.
      \end{aligned}
    \end{equation}
We have then found our upper bound and we can combine (\ref{lower bound proof nuovo 1 thm}) and the result coming from (\ref{first step for upper bound nuovo 1 thm})-(\ref{upperbound proof nuovo primo teorema}) thus obtaining
\begin{equation*}
    \begin{aligned}
    & |W| \sum_{x \in B_{\delta_R}} |u(x)|^p \mu(x) \\
    & \leq \frac{C_2}{\delta^2 R^2}\sum_{x \in B_{R(x_0)}} \mu(x) |u(x)|^p [d(x,x_0) + k]^{-\beta}.
    \end{aligned}
\end{equation*}
Therefore
\begin{equation*}
    |W| \sum_{x \in B_{\delta_R}} |u(x)|^p \mu(x)
    \leq \frac{C_2}{\delta^2 R^2} \|u\|_{l_{\zeta_\beta}^p}^p.
\end{equation*}
Then, letting $R \longrightarrow \infty$,  one has
\begin{equation*}
    |W| \sum_{x \in B_{\delta_R}} |u(x)|^p \mu(x)
    \leq 0.
\end{equation*}
Hence, since both $|W|$ and $\mu(x)$ are positive terms, the above inequality implies
\begin{align*}
    u(x) = 0 \hspace{1.5 cm} \text{for all } x \in G.
\end{align*}
And this concludes the proof.
\end{proof}

\section{Proof of Theorem \ref{teo2}}\setcounter{equation}{0}\label{alt}


We now show the next lemma.
\begin{lemma} {}\label{New Lemma 6.1}
Let $p \geq 1$. Then
\begin{equation*}
\begin{aligned}
    &\Delta \zeta_\beta  - p V \zeta_\beta \leq \bigg[C_0C_1\beta [d(x,x_0) + k]^{-1} - pc_0 [d(x,x_0) + k]^{-\alpha}\bigg] \zeta_\beta  \hspace{0.5 cm} \text{for any } x \in G,
\end{aligned}
\end{equation*}
where $C_1$ and $c_0$ are given in \eqref{C1 DECLARE} and (\ref{ZETA FUNZ PESO NUOVA}) respectively, while and $C_0$ is the bound of the 1-intrinsic metric $d$.
\end{lemma}
\begin{proof}
Recalling the expression of the weighted graph Laplacian one has
\begin{align*}
    \Delta \zeta_\beta(x) = \frac{1}{\mu(x)}\sum_{y \in G} \omega(x,y) \bigg\{ [d(x,x_0) + k]^{-\beta} - [d(y,x_0) + k]^{-\beta} \bigg\},
\end{align*}
thus
\begin{align}\label{last estimate before lemma on difference}
    \Delta \zeta_\beta(x) \leq \frac{1}{\mu(x)}\sum_{y \in G} \omega(x,y) \bigg | [d(x,x_0) + k]^{-\beta} - [d(y,x_0) + k]^{-\beta} \bigg|.
\end{align}
By means of Lemma \ref{Lemma Differenza D(x) e D(y)}, we get
\begin{align*}
     \Delta \zeta_\beta(x) \leq \frac{1}{\mu(x)}\sum_{y \in G} \omega(x,y) \beta [d(x,x_0) + k-s]^{-\beta-1} d(x,y).
\end{align*}
Now, relying on the fact that the metric $d$ is 1-intrinsic with bound $C_0$, one has
\begin{align*}
    \Delta \zeta_\beta(x) \leq C_0 \beta [d(x,x_0) + k-s]^{-\beta-1}.
\end{align*}
Therefore
\begin{align}\label{prima stima interm con pot di ex 6.1}
    \Delta \zeta_\beta(x) - p V \zeta_\beta \leq C_0 \beta [d(x,x_0) + k-s]^{-\beta-1} - p c_0[d(x,x_0) + k]^{-\alpha} [d(x,x_0) + k]^{-\beta}.
\end{align}

In view of \eqref{costante C1 origine} and \eqref{C1 DECLARE}, we get
\begin{align*}
    \Delta \zeta_\beta(x) - p V \zeta_\beta \leq C_0 C_1 \beta [d(x,x_0) + k]^{-\beta-1} - p c_0[d(x,x_0) + k]^{-\alpha} [d(x,x_0) + k]^{-\beta}.
\end{align*}
Grouping the terms one obtains
\begin{align*}
    \Delta \zeta_\beta(x) - p V \zeta_\beta \leq [d(x,x_0) + k]^{-\beta} \bigg[C_0C_1\beta [d(x,x_0) + k]^{-1} - pc_0 [d(x,x_0) + k]^{-\alpha}\bigg].
\end{align*}
\end{proof}

Now, in order to have then all the ingredients for the proof of the theorem, we need first to recall  \cite[Proposition 4.4]{MePu1}.
\begin{proposition} {} \label{second fundamental proposition (4.4)}
Let $u$ be a solution to (\ref{SchrodingerEquation}), let $p \geq 1$ and let $v \in \mathfrak{F}$ be a nonnegative function with finite support. Then
\begin{equation*}
    \sum_{x \in G} |u(x)|^p \{ -\Delta v(x) + p V(x) v(x) \} \mu(x) \leq 0
\end{equation*}
\end{proposition}

\begin{proof}[Proof of Theorem \ref{teo2}]
We can apply Proposition \ref{second fundamental proposition (4.4)} with
\begin{align*}
    v(x) \coloneqq \eta(x) \zeta_\beta(x) \hspace{1 cm} \text{for all } x\in G,
\end{align*}
which is a nonnegative function with finite support.
Now, relying on Lemma \ref{Lemma Prodotto due funz}, we get
\begin{equation}\label{Stima con I3 e I4}
    \sum_{x \in G} |u(x)|^p \{ -\Delta \zeta_\beta(x) + p V(x) \zeta_\beta(x) \} \mu(x) \eta(x) \leq \mathit{J_R ^{(1)}} + \mathit{J_R ^{(2)}},
\end{equation}
where
\begin{align*}
 \mathit{J_R ^{(1)}} \coloneqq \sum_{x \in G} |u(x)|^p \zeta_\beta(x) \Delta \eta(x) \mu(x),  \hspace{0.6 cm}   \mathit{J_R ^{(2)}} \coloneqq  \sum_{x,y \in G} |u(x)|^p (\nabla_{xy}\eta) (\nabla_{xy}\zeta_\beta) \omega(x,y).
\end{align*}
We claim that both $\mathit{J_R ^{(1)}}$ and $\mathit{J_R ^{(2)}}$ tend to $0$ as $R \longrightarrow \infty$.  Let us start by providing an upper bound for  $\mathit{J_R ^{(1)}}$. Obviously
\begin{equation*}
|\mathit{J_R ^{(1)}}| =  \bigg|\sum_{x \in G} |u(x)|^p \zeta_\beta(x) \Delta \eta(x) \mu(x) \bigg|,\qquad |\mathit{J_R ^{(1)}}| \leq \sum_{x \in G} |u(x)|^p \zeta_\beta(x) |\Delta \eta(x)| \mu(x).
\end{equation*}
Then, by using Lemma \ref{Lemma 5.2 articolo parte due}, one has
\begin{align*}
    |\mathit{J_R ^{(1)}}| \leq \frac{C_0}{\delta R} \sum_{x \in G} |u(x)|^p \zeta_\beta(x) \mu(x).
\end{align*}
Now, since $u \in l_{\zeta_\beta}^p$ and since the function $\zeta_\beta$ is actually the weight function, we have that
\begin{align*}
    |\mathit{J_R ^{(1)}}| \leq \frac{C_0}{\delta R} ||u||_{l_{\zeta_\beta}^p}^p \longrightarrow 0,
\end{align*}
as $R$ go to $+\infty$.
Now, we can pass to an estimate for $\mathit{J_R ^{(2)}}$. Clearly,
\begin{equation}\label{First inequality for I4}
|\mathit{J_R ^{(2)}}| \leq  \sum_{x,y \in G} |u(x)|^p |(\nabla_{xy}\eta)| |(\nabla_{xy}\zeta_\beta)| \omega(x,y). \\
\end{equation}
By exploiting Lemma \ref{Lemma Differenza D(x) e D(y)}, we obtain that
\begin{align*}
    |(\nabla_{xy}\zeta_\beta)| \leq \beta d(x,y) [d(x,x_0) + k - s]^{-\beta -1}.
\end{align*}
By the same reasoning used in the proof of Lemma \ref{New Lemma 6.1}, one obtains that
\begin{align*}
    [d(x,x_0) + k-s]^{-\beta-1} \leq C_3 [d(x,x_0) + k]^{-\beta} \hspace{1 cm} \forall x \in G,
\end{align*}
with $\displaystyle
    C_3 \coloneqq \frac{k^\beta}{(k-s)^{\beta+1}}.$
Hence
\begin{align*}
   |(\nabla_{xy}\zeta_\beta)| \leq C_3 \beta d(x,y) [d(x,x_0) + k]^{-\beta }.
\end{align*}
Using the above inequality, (\ref{First inequality for I4}) becomes
\begin{align*}
    |\mathit{J_R ^{(2)}}| \leq  C_3 \beta \sum_{x,y \in G} |u(x)|^p |(\nabla_{xy}\eta)| d(x,y) [d(x,x_0) + k]^{-\beta } \omega(x,y).
\end{align*}
Now, by exploiting Lemma \ref{Lemma 5.2 articolo parte due},
\begin{align*}
     |\mathit{J_R ^{(2)}}| \leq \frac{C_3 \beta}{\delta R}  \sum_{x,y \in G} |u(x)|^p | d(x,y)^2 [d(x,x_0) + k]^{-\beta} \omega(x,y).
\end{align*}
Splitting properly the sum and relying on the definition of $\zeta_\beta$, it is possible to write
\begin{align*}
    |\mathit{J_R ^{(2)}}| \leq \frac{C_3 \beta}{\delta R}  \sum_{x \in G} \zeta_\beta(x) |u(x)|^p \sum_{y \in G} \omega(x,y) d(x,y)^2.
\end{align*}
Thus, relying on the hypothesis for the metric $d$ to be intrinsic, one has
\begin{align*}
    |\mathit{J_R ^{(2)}}| \leq \frac{C_3 \beta}{\delta R} \sum_{x \in G} \zeta_\beta(x) |u(x)|^p \mu(x).
\end{align*}
Since $u \in l_{\zeta_\beta}^p$, one obtains
\begin{align*}
    |\mathit{J_R ^{(2)}}| \leq \frac{C_3 \beta}{\delta R} ||u||_{l_{\zeta_\beta}^p}^p \longrightarrow 0,
\end{align*}
as $R \longrightarrow \infty$. Summing up, we have shown that, letting $R \longrightarrow \infty$, the estimate (\ref{Stima con I3 e I4}) yields
\begin{align*}
    \sum_{x \in G} |u(x)|^p \{ -\Delta \zeta_\beta(x) + p V(x) \zeta_\beta(x) \} \mu(x) \leq 0.
\end{align*}

In view of Lemma \eqref{New Lemma 6.1}, hypothesis \eqref{e2f} and the fact that $\alpha\leq 1$,  the quantity inside the curly brackets is positive. Since also $\mu$ is positive, one necessarily has
\begin{align*}
    u(x) = 0  \hspace{1.5 cm} \text{for all } x \in G.
\end{align*}
This ends the proof of the theorem.
\end{proof}


\end{document}